\documentclass[10pt]{article}

\newcommand{\isd}{\stackrel{d}{=}}

\newcommand{\PP}{\mathbb{P}}
\newcommand{\ZZ}{\mathbb{Z}}

\newcommand{\EE}{\mathbb{E}\,}

\renewcommand{\Box}{\square}

\newcommand{\Ber}{\textrm{Ber}}

\newcommand{\tF}{{\tilde{F}}}
\newcommand{\tf}{{\tilde{f}}}

\usepackage{amsfonts}
\usepackage{amsxtra}
\usepackage{amsmath}
\usepackage{a4}
\usepackage{amssymb}
\usepackage{theorem}
\usepackage{psfrag}
\usepackage{epsfig}
\usepackage{graphicx}
\usepackage{color}

\newtheorem{theorem}{Theorem}

{\theorembodyfont{\rmfamily}}
{\theorembodyfont{\rmfamily}}

\numberwithin{theorem}{section}
\numberwithin{equation}{section}
\numberwithin{figure}{section}

\renewcommand{\Ber}{\ensuremath\textrm{Ber}}
\newcommand{\Geom}{\ensuremath\textrm{Geom}}
\newcommand{\Exp}{\ensuremath\textrm{Exp}}

\allowdisplaybreaks

\begin{document}

\title{Batch queues, 
reversibility and first-passage percolation}
\author{
\textbf{James B.~Martin}
\\
\textit{University of Oxford}}
\date{ }
\maketitle


\begin{abstract}
We consider a model of queues in discrete time, with batch services and arrivals. 
The case where arrival and service batches both have Bernoulli distributions corresponds
to a discrete-time $M/M/1$ queue, and the case where both have geometric 
distributions has also been previously studied. We describe a common extension to
a more general class where the batches are the product of a Bernoulli and a 
geometric, and use reversibility arguments to prove versions of Burke's theorem for these models.
Extensions to models with continuous time or continuous workload are also described. 
As an application, we show how these results can be combined with methods of 
Sepp\"al\"ainen and O'Connell to provide exact solutions for a new class of 
first-passage percolation problems.
\end{abstract}

\section{Introduction}
We consider a model of queues in discrete time, with batch services and arrivals.  
At the beginning of time slot $n$, there are $X_n$ customers in the queue.
A number $A_n$ of customers then arrive, increasing the queue length to $X_n+A_n$. 
After this an amount $S_n$ of service is available, so that $D_n=\min(S_n, X_n+A_n)$ 
customers depart from the queue. Typically, we assume that the sequences
$A_n$ and $S_n$ are random. 

This model has been studied in various contexts (sometimes described as a
\textit{storage} model rather than a queue \cite{DraMaiOCo}).
There is a close correspondence
between this model and another type of queueing model, in which the data $S_n$ represent
inter-arrival times between successive customers and the data $A_n$ represent
service requirements of customers. 

The case where the sequences $S_n$ and $A_n$ are independent and both
consist of i.i.d.\ Bernoulli random variables corresponds to an $M/M/1$ queue 
in discrete time.
The case where the Bernoulli distributions are replaced by geometric distributions
has also been studied previously \cite{BedAzi, DraMaiOCo, OConnell}. 

We generalise these two situations to the case where the distributions are
a product of a Bernoulli and a geometric. 

We show that for appropriate choices of the parameters, the queue is 
reversible in equilibrium, and that various forms of Burke's theorem hold: 
for example, the departure process has the same distribution as the arrival process,
and the queue-length at a given time is independent of the departure process
before that time. The stationary distribution of the queue-length is also 
given by a product of a Bernoulli and a geometric. 

These properties make it easy to describe the stationary behaviour of several
such queues in tandem. As an application, we show how this can be used 
to calculate rates of growth for certain first-passage percolation problems.
This uses techniques developed by Sepp\"al\"ainen \cite{Seppshape} and O'Connell
\cite{OConnell}, 
and extends the class of ``exactly solvable'' 
first-passage percolation models.

In Section \ref{queue} we describe the basic queueing model in more detail.
In Section \ref{BerGeo} we collect some definitions regarding probability distributions,
in particular the distribution obtained by multiplying a Bernoulli
and a geometric. 

In Section \ref{reversibility} we state and prove the reversibility 
results and versions of Burke's theorem, for queues with 
Bernoulli-geometric arrival and service processes. We discuss related 
results concerning systems of queues in tandem. We also note results 
about the stationary distribution of the queue-length in the case where 
only the arrival process is assumed to have this form, and the service 
process consists of any i.i.d.\ sequence. 

As well as extending the Bernoulli and geometric cases, the 
Bernoulli-geometric model has another useful feature, that by taking a 
limit as the parameter of the Bernoulli tends to 0 one can easily arrive 
at various continuous-time models in which arriving customers or offered 
service occur at times corresponding to points of Poisson processes. 
Taking an alternative limit, one can also consider models with 
continuous workload, where arrival and service batches are exponential 
rather than geometric. These extensions are indicated in Section \ref{continuoussection}.

In Section \ref{percolation} we describe the application to 
first-passage percolation models, and give exact expressions for some 
time-constants.

\section{Queueing model with batch services and arrivals}\label{queue}
We describe the main queueing model of the paper. 

The queue is driven by an arrival process $\left(A_n, n\in\ZZ\right)$ and a 
service process $\left(S_n, n\in\ZZ\right)$.

At time-slot $n\in\ZZ$, $A_n$ customers arrive at the queue. 
Then service is available for $S_n$ customers;
if the queue-length is at least $S_n$, then $S_n$ customers are served, while 
if the queue length is less than $S_n$ then all the customers are served. 

Let $X_n$ be the queue length after the service $S_{n-1}$, before the arrival $A_n$.
From the description of the queue above, we have the basic recurrences
\begin{equation}\label{Xrec}
X_{n+1}=\left[ X_n + A_n - S_n\right]_+
\end{equation}
(where $[x]_+$ denotes $\max\{x,0\}$).
Similarly if $Y_n$ is the queue length after the arrival 
$A_n$ and before the service $S_n$, then
$Y_n=X_n+ A_n$, and
\begin{equation}\label{Yrec}
Y_{n+1}=\left[Y_n-S_n\right]_+ + A_{n+1}.
\end{equation}

In this paper we will almost always consider the case where the arrival
and service processes are independent, and $\left(A_n, n\in\ZZ\right)$
and $\left(S_n, n\in\ZZ\right)$ are both i.i.d.\ sequences,
with $\EE A_n<\EE S_n$. In this case (in fact, much more generally) 
we can define the queue-length sequence 
$\left(X_n, n\in\ZZ\right)$ by 
\begin{equation}\label{Xdef}
X_n=\max_{m\leq n}\sum_{r=m}^{n-1}\left(A_r-S_r\right)
\end{equation}
(where a sum from $n$ to $n-1$ is understood to be 0).
This quantity is 
almost surely finite since the common mean of the $S_n$ is larger 
than that of the $A_n$. Then 
the sequence $(X_n)$ satisfies the recurrences (\ref{Xrec}), 
and that $(X_n)$ is a stationary Markov chain. 

Let $D_n$ be the number departing from the queue at the time of the service $S_n$. So
\begin{align*}
D_n&=\min(Y_n, S_n)\\
&=Y_n-X_{n+1}\\
&=X_n+A_n-X_{n+1}\\
&=Y_n+A_{n+1}-Y_{n+1}.
\end{align*}

Let $U_n$ be the unused service at the time of the service $S_n$. So
\begin{align*}
U_n&=S_n-D_n\\
&=[S_n-Y_n]_+.
\end{align*}

\begin{figure}[t]
	\centering
	{\input{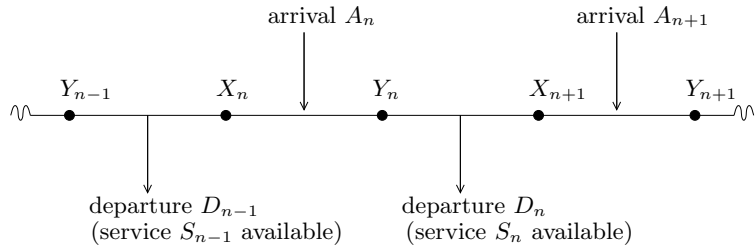}}
	\caption{The evolution of the queue with batch services
and arrivals \label{slotfig}}
\end{figure}

See Figure \ref{slotfig} for a representation of the evolution 
of the queue along with its inputs and outputs. 

We also introduce some further quantities whose interpretation is ostensibly
less natural (but see Section \ref{dual}).
Write $I_n=U_n+A_{n+1}$ for the unused service plus next arrival, 
and $T_n=U_n+A_n$ for the unused service plus previous arrival. 

Note that although we have talked in terms of ``numbers of customers'', 
there is no reason why the variables have to take integer values. 
We will also consider the case where $A_n$ and $S_n$ are non-negative real-valued random 
variables; here one might talk of ``amount of work'', say,  rather than 
``number of customers''.

This model of a discrete-time queue with batch services and arrivals has
been considered in various contexts, for example by 
Bedekar and Azizo\~{g}lu \cite{BedAzi}, 
Ganesh, O'Connell and Prabhakar \cite{GanOCoPra},
Draief, O'Connell and Mairesse \cite{DraMaiOCo}, and
O'Connell \cite{OConnell}.

Models that have been studied include the case where 
arrival batches and service batches are Bernoulli distributed,
and the case where both have a geometric distribution; 
see the next sections for further details.
The main result in this paper is to extend 
the versions of Burke's theorem obtained in these cases 
to the case of a distribution which is the product of a 
Bernoulli and a geometric.

\subsection{Dual queueing model}\label{dual}
Our main batch queueing model is closely related to an alternative
(and in fact more widely-studied) model of a single-server queue 
with first-in-first-out service discipline.

This dual model uses the same variables and recurrences as above, but with different 
interpretations. Now $A_n$ represents the amount of time required by the
$n$th customer for service, and $S_n$ is the interarrival time between customers $n$ and $n+1$.

Let $X_n$ be the waiting time of customer $n$ between his arrival at the queue
and the start of his service. Then $X_n$ obeys the same recurrence as at 
(\ref{Xrec}) above. Now $Y_n$ is the 
total time spent by customer $n$ in the queue
(including service), and $D_n$ is the time spent by customer $n$ 
at the back of the queue. $U_n$ is the idle time of the server 
between departure of customer $n$ and 
arrival of customer $n+1$, $I_n$ is the interdeparture time
between customers $n$ and $n+1$, and $T_n$ is the time between the starts
of service of customers $n$ and $n+1$.

See Draief, Mairesse and O'Connell \cite{DraMaiOCo} for extensive discussions 
of the relations between the two models.

\subsection{Geometric case}\label{geometric}
We assume that the processes $(A_n, n\in\ZZ)$ 
and $(S_n, n\in\ZZ)$ are independent, and
that each process is an i.i.d.\ sequence.

Suppose now that both $A_n$ and $S_n$ have geometric distribution, with 
$\PP(A_n=k)=\alpha(1-\alpha)^{k-1}$ and 
$\PP(S_n=k)=\beta(1-\beta)^{k-1}$ for $k=1,2,\dots$. For stability 
of the process, we require $\beta<\alpha$. 

In this case, the single-server queue of Section \ref{dual} is 
an $M/M/1$ queue in discrete time. A version of Burke's theorem 
was proved for this model by Hsu and Burke in \cite{HsuBurke}.
Among other properties one has that the arrival and departure
processes have the same law, which (in the notation of 
Section \ref{dual}) says that $(S_n, n\in\ZZ)\isd (I_n, n\in\ZZ)$.

On the other hand, Bedekar and Azizo\~glu \cite{BedAzi} 
showed an analogous input-output theorem for the batch-queueing model;
namely, that $(A_n, n\in\ZZ)\isd (D_n, n\in\ZZ)$. 

These results are unified by Draief, Mairesse and O'Connell \cite{DraMaiOCo}, 
who obtain a \textit{joint} Burke's theorem for the two models, 
namely that 
\begin{equation}\label{firstjoint}
((A_n, S_n), n\in\ZZ)\isd ((D_n, I_n), n\in\ZZ).
\end{equation}

The same result also applies if the distributions of $A_n$ and $S_n$ are
exponential rather than geometric.

\subsection{Bernoulli case}
Suppose instead that $A_n$ and $S_n$ both have Bernoulli distribution,
with $p=\PP(A_n=1)=1-\PP(A_n=0)$ and $q=\PP(S_n=1)=1-\PP(S_n=0)$.
For stability of the process, we require $p<q$.

In our main batch-queue model, this means that each arrival batch
is either empty or contains a single customer, and 
at each slot the available service is either 1 or 0. 

In this case, the batch-queue model in fact corresponds to an
$M/M/1$ queue in discrete time (hence to the \textit{dual} queueing model
in the geometric case of Section \ref{geometric}),  since the intervals between 
arrivals of successive customers are geometric 
(with mean $1/p$) and the service time of a customer once 
he reaches the front of the queue is also geometric (with mean $1/q$).
So the statement that $(S_n)\isd(I_n)$ in the geometric case of Section
\ref{geometric} corresponds to the statement that $(A_n)\isd(D_n)$ in 
the Bernoulli case.

Again this can be extended to a sort of joint input-output theorem for the two models
in the Bernoulli case; 
namely one has that 
\begin{equation}\label{secondjoint}
((A_n, S_n), n\in\ZZ)\isd ((D_n, T_n), n\in\ZZ).
\end{equation}
This result is contained in the proof of Theorem 4.1 of 
K\"onig, O'Connell and Roch \cite{KonOCoRoc}. 

Note the difference between (\ref{firstjoint}) and (\ref{secondjoint}). 
In the geometric case, the output theorem for the single-server queue 
involves the interdeparture process $(I_n)$, while in the Bernoulli case
it involves instead the process $(T_n)$ of intervals between \textit{starts} of service.
One can see that $I_n$ is not Bernoulli in the Bernoulli case (indeed, 
one may have $I_n=2$); also, one can see that $T_n$ and $T_{n+1}$ are not 
independent in the geometric case so that (\ref{secondjoint}) certainly does not hold there. 

In Theorem \ref{main}, we'll show that the common part of (\ref{firstjoint})
and (\ref{secondjoint}), namely the 
result that $(A_n)\isd(D_n)$, extends to a class of cases 
where the distributions of $A_n$ and $S_n$ are products 
of a geometric distribution and a Bernoulli distribution.

\section{Bernoulli-geometric distribution}\label{BerGeo}
We first introduce some notation.
$X$ is said to have Bernoulli distribution with parameter $p$ if 
$\PP(X=1)=p$ and $\PP(X=0)=1-p$.

We say that $X$ has $\Geom^+$ distribution with parameter $\alpha\in(0,1)$ if 
\[
\PP(X=k)=\alpha(1-\alpha)^{k-1}
\]
for $k\geq 1$. If $X$ has $\Geom^+(\alpha)$ distribution then 
$X-1$ is said to have $\Geom^0(\alpha)$ distribution.

Now we define a Bernoulli-geometric distribution, with parameters $p$ 
and $\alpha$.
A random variable with this distribution has the distribution of the product of 
two independent random variables, one with $\Ber(p)$ distribution and the other with 
$\Geom^+(\alpha)$ distribution. That is, 
$A\sim\Ber(p)\Geom(\alpha)$ if 
\[
\PP(A=k)=\begin{cases}
1-p,& k=0\\
p\alpha(1-\alpha)^{k-1},& k\geq 1.\end{cases}
\]
We have $\EE A=p/\alpha$, and the probability generating function of $A$ is given by  
\[
\EE(z^A)=\frac{(1-p)-(1-p-\alpha)z}{1-(1-\alpha)z}.
\]
The distribution of $A$, conditioned on being non-zero, is simply $\Geom^+(\alpha)$.

In passing, note that 
such a random variable 
$A$ may also be represented as a geometric number of independent
geometrics, in the case $p<\alpha$.
Namely, let $V\sim\Geom^0(p)$, and 
let $W_i$ be i.i.d.\ $\Geom^+(\gamma)$, 
where $\gamma=(\alpha-p)/(1-p)$,
and independent of $V$.  
Then define $R=W_1+W_2+\dots+W_V$. One has $R\stackrel{d}{=}A$.

Alternatively, let $V\sim\Geom^0(r)$ where $r=(1-p)p/(1-p\alpha)$,
and let $W_i$ be i.i.d.\ $\Geom^0(\gamma)$, where $\gamma$ is as above,
and independent of $V$. As before, define $R=W_1+\dots+W_V$; again one gets $R\stackrel{d}{=}A$.

\section{Bernoulli-geometric queues}\label{reversibility}
We consider a queue with arrival process $A_n, n\in\ZZ$ 
where $A_n$ are i.i.d. Ber($p$)Geom($\alpha$),
and a service process $S_n, n\in\ZZ$ independent of the arrival process
and with $S_n$ i.i.d. Ber($q$)Geom($\beta$).

For stability we clearly need $\EE S_n>\EE A_n$, i.e.\ $p\beta<q\alpha$.

In general, the queue-length process is not reversible. For example,
if $\alpha$ is small and $\beta$ is large, then one tends to see 
infrequent but large arrival batches, and frequent but small departure batches. 

However, under a further condition, which in effect
reduces the number of parameters from 4 to 3,  we show reversibility and 
various other related properties. 
The relevant condition is that
\begin{equation}\label{1param}
\frac{\alpha}{1-\alpha}\frac{p}{1-p}=\frac{\beta}{1-\beta}\frac{q}{1-q}.
\end{equation}
Note that combined with the stability condition $p\beta<q\alpha$, 
condition (\ref{1param}) implies that $\beta<\alpha$ and $p<q$.

For a given service distribution, condition (\ref{1param}) 
gives one degree of freedom
for the arrival distribution. In particular, suppose $q$ and $\beta$ are fixed,
giving an overall 
service intensity of $q/\beta$. Then for any $\lambda<q/\beta$,
there exists a unique pair $p$, $\alpha$ satisfying (\ref{1param}) 
and such that the arrival intensity $p/\alpha$ is equal to $\lambda$.
We further discuss the relevance of (\ref{1param}) to the reversibility of the queueing process 
in Section \ref{1paramdiscusssection}.

\begin{theorem}\label{main}
Suppose that $p\beta<q\alpha$, and (\ref{1param}) holds.
\begin{itemize}
\item[(i)] The queue-length processes $(X_n)$ and $(Y_n)$ are reversible;
moreover, they are jointly reversible in the sense that
\begin{equation}
\label{reverse}
\left(\dots,X_{-1},Y_{-1},X_0,Y_0,X_1,Y_1,\dots\right)
\stackrel{d}{=}
(\dots,X_2,Y_1,X_1,Y_0,X_0,Y_{-1},\dots).
\end{equation}
\item[(ii)]The departure process $\left(D_n, n\in\ZZ\right)$ has the same
law as the arrival process $\left(A_n, n\in\ZZ\right)$.
\item[(iii)]The stationary distributions of the queue-length processes 
(before and after service) are given by
\[
X_n\sim\Ber\left(c\right)
\Geom\left(\gamma\right),
\]
where 
\[c=\frac{\beta}{1-\beta}\frac{1-\alpha}{\alpha},\,\,\,\gamma=\frac{\alpha-\beta}{1-\beta},\]
and 
\[
Y_n\sim \Ber(p+c-pc)\Geom\left(\gamma\right).
\]
\item[(iv)] For all $n$, the queue length $X_n$ at time $n$ is independent of the 
process of departures $\left(D_i, i<n\right)$ before time $n$. 
\end{itemize}
\end{theorem}

From (iii), one has 
\[
\EE X=\frac{\beta(1-\alpha)}{\alpha(\alpha-\beta)},\,\,\,\,
\EE Y=\EE X+\frac p\alpha.
\]

\noindent\textit{Proof of Theorem \ref{main}}.
To show the reversibility in (i) and the stationary distribution for $X_n$ in (iii),
it is enough to show that for 
all $k,m,r$, with $m\geq k$ and $m\geq r$,
\begin{equation}\label{tocheck}
\pi(k)\PP(Y_0=m|X_0=k)\PP(X_1=r|Y_0=m)
=
\pi(r)\PP(Y_0=m|X_0=r)\PP(X_1=k|Y_0=m),
\end{equation}
where $\pi$ is the probability mass function for the distribution of $X_n$ given in (iii).

If $k=r$, this is obvious. The cases $k<r$ and $k>r$ are symmetric, and one only needs to 
check one, say $k<r$. 

Note that 
\[
\PP(Y_0=m|X_0=k)=\PP(A=m-k),
\]
while
\[
\PP(X_1=r|Y_0=m)=\begin{cases}
\PP(S\geq m-r), &r=0\\
\PP(S=m-r),&r>0
\end{cases}.
\]
Further one has by assumption that 
\begin{align*}
\pi(x)=
\begin{cases}
1-c, &x=0\\
c\gamma(1-\gamma)^{x-1}, &x\geq 1
\end{cases}&,\,\,\,\,\,\,\,\,
\PP(A=x)=\begin{cases}
1-p, &x=0\\
p\alpha(1-\alpha)^{x-1}, &x\geq 1
\end{cases},
\\
\PP(S=x)=\begin{cases}
1-q, &x=0\\
q\beta(1-\beta)^{x-1}, &x\geq 1
\end{cases}&,\,\,\,\,\,\,\,\,
\PP(S\geq x)=\begin{cases}
1, &x=0\\
q(1-\beta)^{x-1}, &x\geq 1
\end{cases}.
\end{align*}

Now one can for example divide into four further cases: (1) $k>0, m=r$; 
(2) $k>0, m>r$; (3) $k=0, m=r$; (4) $k=0, m>r$, and check (\ref{tocheck}) directly
in each case. 

For example, in case (1), we can use the forms of the probability distributions above to 
give 
\begin{multline*}
\pi(k)\PP(Y_0=m|X_0=k)\PP(X_1=m|Y_0=m)\\
\begin{aligned}
&=\pi(k)\PP(A=m-k)\PP(S=0)\\
&=c\gamma(1-\gamma)^{m-1}(1-p)q\beta(1-\beta)^{m-k-1}
\left[\frac{p}{1-p}\frac{\alpha}{1-\alpha}\frac{1-q}{q}\frac{1-\beta}\beta\right]
\left[\frac{1-\alpha}{(1-\beta)(1-\gamma)}\right]^{m-k}\\
&=c\gamma(1-\gamma)^{m-1}(1-p)q\beta(1-\beta)^{m-k-1}\\
&=\pi(m)\PP(A=0)\PP(S=m-k)\\
&=\pi(m)\PP(Y_0=m|X_0=m)\PP(X_1=k|Y_0=m).
\end{aligned}
\end{multline*}
Two lines from the end, we used condition (\ref{1param}) and the fact that
$\frac{(1-\beta)(1-\gamma)}{1-\alpha}=1$ (which follows from the definition of $\gamma$).

The other three cases follow similarly and we omit the details.

The stationary distribution for $Y_n$ in (iii) follows from 
the distribution of $X_n$ and the fact that $Y_n=X_n+A_n$ with $X_n$ and $A_n$ 
independent (for example, multiply the generating functions).

Given the reversibility, the properties in (ii) and (iv) can be deduced
using the same arguments that have been used to establish Burke's theorem in 
various settings (originally by Reich for the continuous-time $M/M/1$ queue
\cite{Reich}).

First note that $D_n=Y_n-X_{n+1}$ and $A_n=Y_n-X_n$. Hence the reversibility property
in (i) implies that $(D_n)_{n\in\ZZ}$ and $(A_{-n})_{n\in\ZZ}$ have the same distribution.
But $A_n$ are i.i.d.\ random variables, so $(A_{-n})_{n\in\ZZ}$ and $(A_n)_{n\in\ZZ}$
in turn have the same distribution, giving (ii). 

Now note that, by (\ref{Xdef}) for example, we can write $X_n$ as a function
of $(A_i)_{i<n}$ and $(S_i)_{i<n}$. Hence $X_n$ is independent of $(A_i)_{i\geq n}$. 
But, using (i) again, the collection
\[
\left(X_n; A_n, A_{n+1}, A_{n+2},\dots) 
= (X_n; Y_n-X_n, Y_{n+1}-X_{n+1}, Y_{n+2}-X_{n+2},\dots\right)
\]
has the same distribution as the collection
\[
\left(X_n; D_{n-1}, D_{n-2}, D_{n-3}, \dots) 
= (X_n; Y_{n-1}-X_n, Y_{n-2}-X_{n-1}, Y_{n-3}-X_{n+2}, \dots\right).
\]
Thus indeed $X_n$ is independent of $(D_i)_{i<n}$ as required for (iv).
$\hfill\Box$

\subsection{Fixed points}
Consider a queueing server defined by a given distribution of the service process.
We may ask for distributions of the arrival process with the following property:
when such an arrival process is fed into the queue (independently of 
the service process), the resulting departure process has the same law as the arrival process.
Such a distribution of the arrival process is called a fixed point for the 
given service process.

Starting from Burke's theorem for a continuous-time $M/M/1$ queue,
questions concerning fixed points have been extensively studied 
in the context of single-server queues such as the model of Section \ref{dual} --
see for example \cite{MaiPra} and references therein. 
They have also been considered, although less often, in the context
of the model of batch arrivals and services in discrete time considered here --
see for example \cite{GanOCoPra}.

Part (ii) of Theorem \ref{main} is such a \textit{fixed point} result.
Fix a service process of $\Ber\Geom$ type, specified by the parameters $q$ and $\beta$.
Let $\mu=q/\beta=\EE S_n$ be the service intensity. 
Then for each $\lambda<\mu$, there exists an arrival process which is a fixed point of the queue,
and which has arrival intensity $\EE A_n=\lambda$ 
(simply choose $p$ and $\alpha$ to satisfy $p/\alpha=\lambda$ along with
condition (\ref{1param}) -- there is a unique way to do this).

In fact, Theorem 5 of \cite{GanOCoPra} implies that this gives
the \textit{unique} fixed-point arrival process which is ergodic with
arrival intensity $\lambda$. Furthermore, this fixed point is
\textit{attractive}; loosely, this means that if any ergodic arrival
process is fed into a tandem of queueing servers with 
this service distribution, the distribution of the resulting output
process converges to the fixed point as the length of the tandem grows. 
See \cite{GanOCoPra} for precise definitions.

We note in passing that if (\ref{1param}) is satisfied, then
the distribution $\Ber(p)\Geom(\alpha)$ has minimal relative entropy
with respect to the distribution $\Ber(q)\Geom(\beta)$, 
out of all distributions with mean $p/\alpha$. See the introduction
of \cite{GanOCoPra} for related discussions.

\subsection{Tandems}\label{tandem}
Using Theorem \ref{main} we can also describe the behaviour 
of systems of queues in tandem. 

Consider a system of $R$ queues in tandem. Each queue has an independent
service process, $S_n^{(r)}$ for the $r$th queue, and each of these processes
is a collection of i.i.d.\ $\Ber(q)\Geom(\beta)$ random variables.

The first queue has an arrival process $A_n^{(1)}$, which is independent of 
the service processes and is a collection of i.i.d.\ $\Ber(p)\Geom(\alpha)$
random variables. We assume that $p$, $q$, $\alpha$, $\beta$ satisfy $(\ref{1param})$
and the stability condition $p\beta<q\alpha$.

Now, recursively, let the arrival process to the $r$th queue be given by the 
departure process from the $(r-1)$st queue, for $r=2,3,\dots,R$; 
that is, $A_n^{(r)}=D_n^{(r-1)}$. Thus a customer departing from queue $r-1$
moves immediately (within the same time-slot) to the next queue.

Using Theorem \ref{main} and well-known methods, we obtain the following result:
\begin{theorem}\label{tandemtheorem}
\begin{itemize}
\item[(i)]
All the departure processes $D^{(r)}$ have the same distribution, which is also
the distribution of $A^{(1)}$.
\item[(ii)]
The vector $X^{(1)}_n, X^{(2)}_n, \dots, X^{(R)}_n$ of queue-lengths
of the $R$ queues at a fixed time $n$ is a collection of i.i.d.\ random variables, whose
common distribution is that
given in Theorem \ref{main}(iii). 
\end{itemize}
\end{theorem}

\noindent\textit{Proof:}
Part (i) is obtained by applying Theorem \ref{main}(ii) repeatedly.
The argument to obtain the product form result in (ii) from Theorem \ref{main}(iv)
is exactly the same as for the familiar case of $M/M/1$ queues in tandem --
see for example Section 2.2 of \cite{Kellybook}.$\hfill\Box$

The result also extends easily to cases where the parameters
of the service process vary between queues: $S_n^{(r)}\sim\Ber(q_r)\Geom(\beta_r)$.
However, all the pairs $(q_r,\beta_r)$ must still belong to the same one-parameter 
family satisfying (\ref{1param}).

We could also consider a vector of ``queue-lengths before service'', $Y^{(r)}_n$. 
Observe that the way we have defined the system of queues in tandem, a customer
may be present after arrival and before service in several different queues at the same
time-slot (since a departure from queue $r-1$ at time $n$ arrives at queue $r$ at the 
same time $n$). So in this case, the corresponding result is in fact that 
$Y_n^{(1)}, Y_{n-1}^{(2)},\dots,Y_{n-R+1}^{(R)}$ form an i.i.d.\ sequence. 
This may also be proved by similar methods.

\subsection{General service-batch distributions}
In order for the reversibility properties in Theorem \ref{main}(i),(ii),(iv) to hold,
one needs the condition (\ref{1param}) relating the parameters of the distributions
of arrival and service distributions. However, one may wonder whether this is necessary
to have a result like Theorem \ref{main}(iii) on the stationary distribution of the
queue-length.

In fact, such a property holds in a much more general case. We will still assume
the same sort of arrival process, but now we will consider \textit{any} service
process which has i.i.d.\ entries which are non-negative integers.
\begin{theorem}\label{generaltheorem}
Suppose $(S_n)$ and $(A_n)$ are both i.i.d.\ sequences taking
non-negative integer values, and independent of each other,
with $\EE S_n>\EE A_n$. 
\begin{itemize}
\item[(a)] If $A_n$ has a $\Ber()\Geom()$ distribution, then both $X_n$ and $Y_n$ have
$\Ber()\Geom()$ distributions.
\item[(b)] If $A_n$ takes values 0 and 1 only, then $X_n$ has a $\Geom^0$ distribution.
\item[(c)] If $A_n$ has a $\Geom^0$ (respectively $\Geom^+$) distribution,
then also $Y_n$ has a $\Geom^0$ (respectively $\Geom^+$) distribution.
\end{itemize}
\end{theorem}
Parts (b) and (c) are both special cases of part (a). Part (c) was already observed
in Proposition 12 of \cite{BedAzi}, and results like part (b) for $M/GI/1$ queues
are certainly well known. 

\medskip
\noindent\textit{Proof:}
We will use a representation of the queue length as the future maximum of a random walk. 
Arguments of this kind are rather classical; for example, see the book of Tak\`acs 
\cite{Takacsbook} for many 
examples (often in the context of queueing theory).

Using (\ref{Xdef}), we can write
\[
X_0=\max_{m\geq 0}\sum_{r=1}^m\left(-S_{-r}+A_{-r}\right), 
\,\,\,\,\,\,\,\,\,
Y_0=A_0+\max_{m\geq 0}\sum_{r=1}^m\left(-S_{-r}+A_{-r}\right).
\]

Consider for example part (b). $X_0$ is the maximum of the walk
\[
0, -S_{-1}+A_{-1}, -S_{-1}+A_{-1}-S_{-2}+A_{-2}, \dots.
\]
The maximum is almost surely finite since the walk has negative drift.
Any positive step of this walk has size exactly 1 (since the $A_n$ are 0 or 1 and the 
$S_n$ are non-negative). Hence any new maximum must be precisely 1 greater than the
previous maximum. 
We may treat these maxima as renewal times; the future evolution of 
the walk (relative to its current position) has the same distribution as the original
walk. Note that $X_0\geq k$ iff this walk reaches level $k$ at some point. 
Hence we get 
\[
\PP(X_0\geq k+1|X_0\geq k)=\PP(X_0\geq 1|X_0\geq 0)=\PP(X_0\geq 1)
\]
for all $k$. Thus $X_0$ indeed has a geometric distribution as desired (and by stationary
the same is true of $X_n$ for all $n$).

For parts (a) we generalise this argument slightly. 
Now the arrival batches may have size greater than 1, but
we will regard such a batch as a sequence of individual steps up, each of size 1.
Since the arrival batches are $\Ber\Geom$, we have the memoryless property for $A_n$:
the distribution of $A_n-k$, conditional on $A_n\geq k$, is the same for all $k\geq 1$.
We regard the walk as a sequence of steps up of size 1, separated perhaps by some number of 
jumps down (which may be 0). By the memoryless property, the jumps down which separate 
each pair of steps up form an i.i.d. sequence, and so we have a 
renewal property after each step up, and thus in particular after each new maximum. 
Thus $\PP(X_0\geq k+1|X_0\geq k)$ doesn't depend on $k$, for $k\geq 1$, and 
indeed $X_0$ has a $\Ber\Geom$ distribution. (The difference from the previous
paragraph is that the distribution for $k=0$ may be different; hence we get a
$\Ber\Geom$ distribution in general rather than the $\Geom$ distribution we had
in the specific case above). The same argument applies also to $Y_0$.

In the special case (c), 
the arrival batches are $\Geom^+$ so that $Y_0$ must always be strictly positive;
hence in fact $Y_0$ is itself $\Geom^+$. If the arrival batches
are $\Geom^0$, we could add 1 to every arrival batch and to every service batch 
to arrive at the previous case; one obtains that $Y_0+1$ is $\Geom^+$ and hence that
$Y_0$ has a $\Geom^0$ distribution. $\hfill\Box$

\subsection{Continuous models}
\label{continuoussection}
All the results above have equivalent versions in the case where 
geometric distributions are replaced by exponential distributions. 

$X$ is said to have $\Ber(p)\Exp(\alpha)$ distribution if
\[
\PP(X\geq x)=\begin{cases}
1,& x=0\\
pe^{-\alpha x},& x>0.\end{cases}
\]
$X$ may be represented as the product of a Bernoulli random variable and an exponential random 
variable, or as the sum of a geometric number of i.i.d.\ exponential random variables.

Now one could write versions of Theorems \ref{main} and \ref{tandemtheorem} where the 
$\Ber()\Geom()$ distributions are replaced by $\Ber()\Exp()$ distributions. These 
continuous versions can be
proved directly, or derived from the discrete versions 
by taking an appropriate limit under which the parameter of the 
geometric distribution tends to 0. 
The analogous condition in place of (\ref{1param}) is that $\alpha p/(1-p)=\beta q/(1-q)$.
Theorem \ref{generaltheorem} extends similarly (and 
we no longer need to assume that the services $S_n$ take integer values). 

An alternative extension is to systems in continuous time rather than discrete time. 
Write the Bernoulli parameters as $p=\epsilon\lambda$ and $q=\epsilon\mu$; now let
$\epsilon\to0$ and rescale time by $\epsilon$. This produces a model in which 
arrival and service batches occur at the times of independent Poisson processes, 
with parameters $\lambda$ and $\mu$ respectively. Again, analogous versions of all our 
main results can be straightforwardly formulated. Note that in this case, the distinction
between queue-length after service and queue-length before service disappears. 

These continuous limits, both in time and in workload, are also exploited in the next section
in the context of first-passage percolation models. 

\subsection{Role of condition (\ref{1param})}
\label{1paramdiscusssection}
It may be helpful to give some brief comments about the relevance of condition 
(\ref{1param}) to the reversibility of the queue-length process as in Theorem \ref{main}.

Consider a ``busy period'' of the process, i.e.\ an excursion away from 0. 
Suppose that the sequence of arrival batches in the busy period is $a_1, a_2, \dots, a_n$
and the sequence of departures is $d_1, d_2, \dots, d_n$. We have $a_1+\dots+a_n=d_1+\dots+d_n$,
and that $a_1+\dots+a_m>d_1+\dots+d_m$ for $m=1,2,\dots, n-1$. Note that this also gives
in particular that  $a_1>0$ 
and $d_n>0$. 

For reversibility of the queue, we need the likelihood of such an excursion to be invariant
under time-reversal (in which the role of arrivals and departures is exchanged).
This likelihood is given by 
\begin{gather*}
\PP\Big( 
X_0=0, A_1=a_1, \dots, A_n=a_n, D_1=d_1, \dots, D_n=d_n\Big),
\\
\intertext{which may also be written as}
\PP\Big(X_0=0, A_1=a_1, \dots, A_n=a_n, S_1=d_1, \dots, S_{n-1}=d_{n-1}, S_n\geq d_n\Big).
\end{gather*}
Note that when we translate from the variables $D_i$ to the variables $S_i$, 
the last equality becomes an inequality, since it is at this point 
that the queue-length returns to 0 and 
so some part of the final service $S_n$ may be unused.

First consider the case where the arrivals are i.i.d.\ $\Geom^0(\alpha)$ 
and the services are i.i.d.\ $\Geom^0(\beta)$. Then the likelihood above 
is equal to 
\[
\alpha^n(1-\alpha)^{a_1+\dots+a_n}\beta^{n-1}(1-\beta)^{d_1+\dots+d_n}.
\]
Since the sum of the $a_i$ is equal to the sum of the $d_i$, this 
is invariant under the exchange $(a_1,\dots,a_n)\leftrightarrow(d_n,\dots,d_1)$ 
as required.

Now consider the case where $A_i\sim\Ber(p)\Geom(\alpha)$ and $S_i\sim\Ber(q)\Geom(\beta)$. 
Now the likelihood above becomes
\[
p^n\alpha^n(1-\alpha)^{a_1+\dots+a_n}
\left(\frac{1-p}p\frac{1-\alpha}\alpha\right)^{\# i:a_i=0}
q^n\beta^{n-1}(1-\beta)^{d_1+\dots+d_n}
\left(\frac{1-q}q\frac{1-\beta}\beta\right)^{\# i:d_i=0}.
\]
Since in general the number of $a_i$ which are zero may be different from the number 
of $d_i$ which are zero, this likelihood is invariant under the exchange 
$(a_1,\dots,a_n)\leftrightarrow(d_n,\dots,d_1)$ only if condition (\ref{1param}) holds.
Hence (\ref{1param}) is also necessary for reversibility of the queue-length process.

By decomposing the queue-length process into its excursions and arguing in this way,
we could in fact arrive at a proof of Theorem \ref{main} which avoided many of the calculations
needed in the proof given above (at the expense of complicating the structure 
of the proof a little). 

The following property is also related to the discussion above. 
If $A\sim\Geom^0(\alpha)$, then the ratio
$\PP(A=k+1)/\PP(A=k)$ is the same for all $k\geq0$, namely $(1-\alpha)$. If instead
$A\sim\Ber(p)\Geom(\alpha)$, then the same is true \textit{except} for $k=0$; 
in the case $k=0$, the ratio is multiplied by a further factor 
$\dfrac{\alpha}{1-\alpha}\dfrac{p}{1-p}$. Condition (\ref{1param}) says that 
this ``adjustment factor'' is the same for the distribution of arrivals as it is
for the distribution of services. 

\section{First-passage percolation models}\label{percolation}
In this section we consider various directed first-passage percolation models,
and use Theorems \ref{main} and \ref{tandemtheorem} to calculate
exact values for time constants.

\begin{figure}[ht]
\label{pathfig}
\begin{center}
\includegraphics[width=0.3\linewidth]{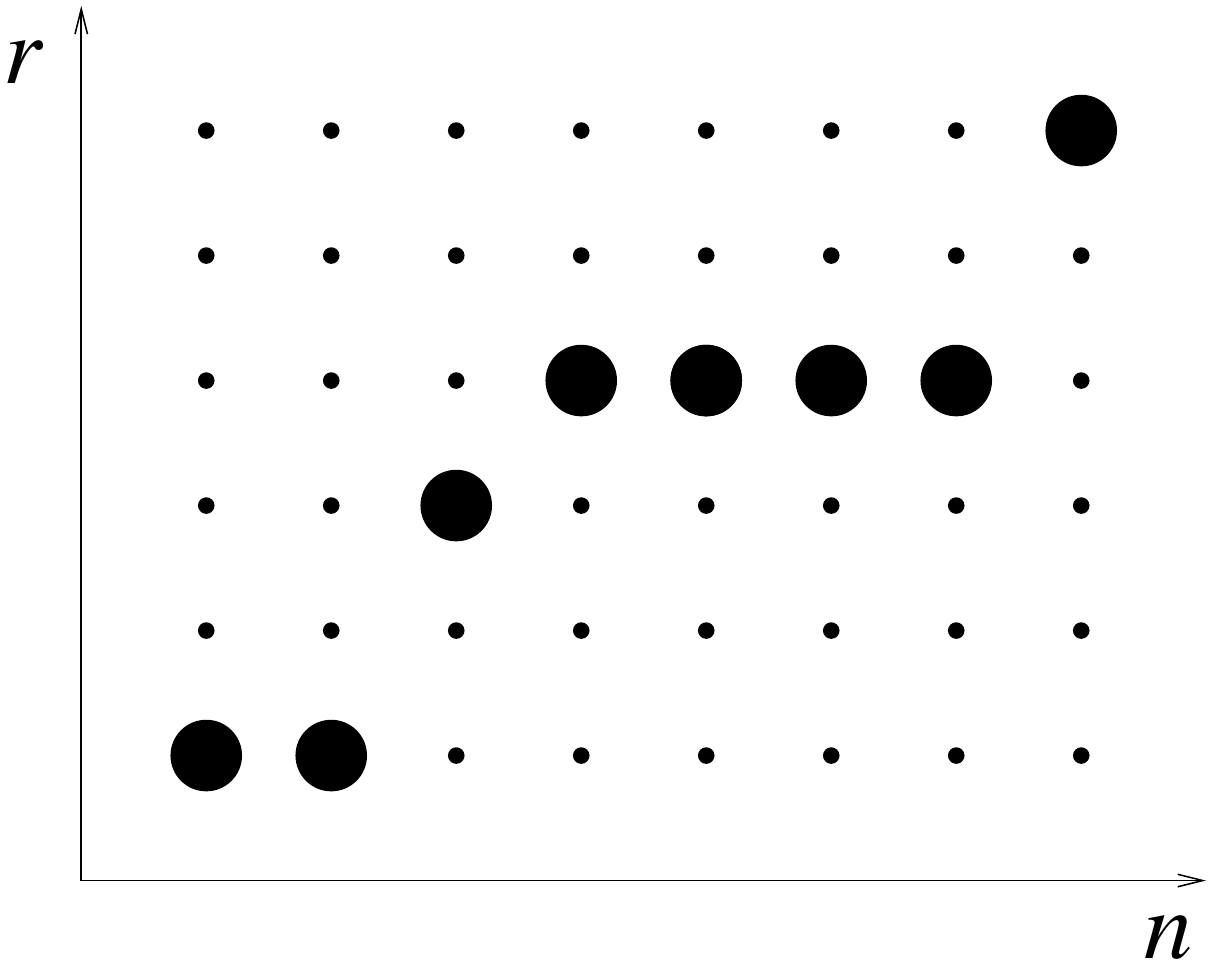}
\caption{An example of a directed path from $(1,1)$ to $(8,6)$.}
\end{center}
\end{figure}

For $(i,j)\leq(k,l)\in\ZZ^2$, denote by $\Pi\left((i,j),(k,l)\right)$ 
the set of ``directed paths''
$(i,r_1), (i+1, r_2), \dots, (k, r_{k-i+1})$, where
$l\leq r_1\leq r_2\leq \dots\leq r_k-i+1\leq l$. 
These paths are strictly increasing in the first coordinate, and 
weakly increasing in the second coordinate. See Figure \ref{pathfig} 
for an example.

Let $S^r_n$ be a collection of i.i.d.\ random variables, and for each 
path $\gamma\in\Pi\left((i,j),(k,l)\right)$,
define the weight of the path $\gamma$ by $S(\gamma)=\sum_{(n,r)\in\gamma}S^r_n$.
Now we define the first-passage time from $(i,j)$ to $(k,l)$ as the minimum
weight over all paths from $(i,j)$ to $(k,l)$:
\begin{equation}\label{Fdef}
F\left((i,j),(k,l)\right)=\min_{\gamma\in\Pi\left((i,j),(k,l)\right)}S(\gamma).
\end{equation}

This model has various possible alternative presentations, for example
in \cite{Seppshape} as a first-passage percolation model with weights on the edges, 
where weights on vertical edges are all equal to a constant and weights on horizontal 
edges are i.i.d. 

By Kingman's subadditive ergodic theorem, for any $x>0$, there exists
a constant $f(x)$ such that
\begin{equation}\label{fdef}
\frac1N F\left((0,0),(\lfloor xN\rfloor,N)\right)\to f(x)\
\end{equation}
almost surely. 

In \cite{Seppshape}, it is shown that if the random variables $S^r_n$ 
have $\Ber(q)$ distribution, then
\begin{equation}\label{Bercase}
f(x)=\begin{cases} 
0,&x\leq(1-q)/q,\\
\left(\sqrt{qx}-\sqrt{1-q}\right)^2,&x>(1-q)/q.
\end{cases}
\end{equation}
In \cite{OConnell}, related methods are used to show that
if the $S^r_n$ have $\Geom^0(\beta)$ distribution, then 
\begin{equation}\label{Geomcase}
f(x)=\begin{cases}
0,&x\leq\beta/(1-\beta),\\
\frac1\beta\left(\sqrt{1-\beta}\sqrt{1+x}-1\right)^2,&x>\beta/(1-\beta),
\end{cases}
\end{equation}
while if the $S^r_n$ have $\Exp(1)$ distribution then
\begin{equation}\label{Expcase}
f(x)=\left(\sqrt{1+x}-1\right)^2.
\end{equation}

The method of \cite{OConnell} exploits a link between
the percolation model and the system of batch queues in tandem of 
the sort described in Section \ref{tandem}. Using this method, we can 
extend the results above to the case of a $\Ber()\Geom()$ distribution.
As suggested by the notation, 
the weights in the percolation problem correspond to service times in 
the queueing model. 

Fix $q$ and $\beta$ and assume that all the $S^r_n$ are i.i.d.\ 
$\Ber(q)\Geom(\beta)$. In this case the service rate at each queue is
$\mu=q/\beta$. For any $\lambda<\mu$, we can choose $p$ and $\alpha$
satisfying (\ref{1param}) and such that $p/\alpha=\lambda$. 
If the arrival batches are i.i.d.\ $\Ber(p)\Geom(\alpha)$ then 
the arrival rate is $\lambda$.

(Rather than choosing $\lambda<\mu$ directly, we may equivalently choose
$\alpha>\beta$ or choose $p<q$. To express one of the three variables $\alpha,p,\lambda$ 
in terms of another, we have, as well as (\ref{1param}), that
\begin{gather}\label{lambdaeq}
\lambda=\frac p\alpha=p\left[\frac{p}{1-p}\frac{1-q}q\frac{1-\beta}\beta+1\right],\\
\lambda=\frac{(1-\alpha)\beta q}{\alpha^2(1-\beta-q)+\alpha\beta q}.
\end{gather}
Recall that $\beta$ and $q$ are to be regarded as fixed throughout).

Now the results of Theorems \ref{main} and \ref{tandemtheorem}
apply. In particular, the queue lengths $X^1_0,\dots,X^R_0$ are i.i.d.,
and their expectation is given by
\begin{equation}\label{heq1}
\EE X^r_0=\frac{\beta}{\alpha-\beta}\frac{1-\alpha}{\alpha}.
\end{equation}
We may regard this expectation as a function of the arrival rate 
$\lambda$ and write it as $h(\lambda)$ to emphasise this. We also have
\begin{equation}\label{heq2}
h(\lambda)=\frac p{(1-p)} \frac{(1-q)}q\left[\frac{p(1-q)}{\beta(q-p)}+1\right].
\end{equation}

Now, in Section 3 of \cite{OConnell}, it is shown, by working recursively
from representations such as (\ref{Xdef}), that 
\[
\sum_{r=1}^R X^r_0=
\sup_{m\leq 0}
\left\{
\sum_{r=m}^{-1} A_r - F\left((m,1),(-1,R)\right)\right\},
\]
and that by passing to the limit
$R\to\infty$ and taking a Legendre transform, one obtains
\begin{equation}\label{feq}
f(x)=\sup_{0<\lambda<\mu}\left\{\lambda x - h(\lambda)\right\}.
\end{equation}
We apply this result to various cases in turn. 

\medskip

\noindent\textbf{Example 1.}
The weights $S^r_n$ have $\Ber(q)\Geom(\beta)$ distribution.
This is essentially the most general case of the ones we treat:
all the others will be derived by taking some sort of limit from this case.

Using (\ref{feq}) and plugging in (\ref{lambdaeq})-(\ref{heq2}), we get
\begin{gather}\label{generalcase}
f(x)=
\frac1{\beta q}
\sup_{p\in(0,q)}
\left\{
\frac{p\big[p(1-q)+(q-p)\beta\big]}{(1-p)}\left[x-\frac{1-q}{q-p}\right]
\right\}\\
\intertext{or alternatively}
f(x)=
\beta
\sup_{\alpha\in(\beta,1)}
\frac{1-\alpha}{\alpha}
\left[
\frac{qx}{\alpha(1-\beta-q)+\beta q}-\frac1{\alpha-\beta}
\right].
\end{gather}
In principle one could solve a quartic equation to put these expressions into
closed form, but the result is unlikely to be pretty. However, it is 
straightforward to show that 
$f(x)=0$ if $x<(1-q)/q$ and $f(x)>0$ otherwise.

\medskip

\noindent\textbf{Example 2.}
Now consider the case where the common distribution of the weights
is $\Ber(q)\Exp(1)$. To obtain this case 
we can let $\beta\to0$ in the previous example and multiply by $\beta$.
We obtain
\[
f(x)=
\sup_{0<r<1} r^2\left[\frac{qx}{1-q+rq}-\frac1{1-r}\right].
\]

\medskip

\noindent\textbf{Example 3.}
Now we consider a model where space becomes continous in one direction.
Consider taking the Bernoulli parameter $q$ to 0, and the space parameter $x$ 
to infinity. We arrive at a model where the first space parameter $i$ is 
replaced by a continuous parameter. For each $r$, we have a ``service process''
$S^r(t)$. This process is a jump process; events occur at times of a Poisson process
of rate 1, say, and to each event is associated a weight, which is the
amount by which the process $S^r$ jumps up at the time of the event. These weights
are independent and each has $\Geom(\beta)$ distribution. In the queueing model,
the weight occurring at time $t$ 
in the process $S^r$ corresponds to the amount of service available at queue $r$ at time $t$. 

The first-passage time can now be defined by 
\begin{equation}\label{Ftildedef}
\tF\left((s,j),(t,l)\right)=
\inf_{s=u_j<u_{j+1}<\dots<u_{l}<u_{l+1}=t}
\sum_{r=j}^l
\left[S^r(u_{j+1})-S^r(u_{j})\right],
\end{equation}
and the time constants by 
\begin{equation}\label{ftildedef}
\tf(y)=
\lim_{N\to\infty}
\frac1N\tF\Big(\big(0,0\big),\big(\lfloor Ny\rfloor,N\big)\Big),
\end{equation}
which, as at (\ref{fdef}), are a.s.\ constant for each $y$.

To obtain this case from (\ref{generalcase}), we let $q\to0$ and set $x=y/q$.
We obtain
\[
\tf(y)=
\beta\sup_{\alpha\in(\beta,1)}\frac{1-\alpha}\alpha
\left[\frac{y}{\alpha(1-\beta)}-\frac1{\alpha-\beta}\right].
\]

\medskip

\noindent\textbf{Example 4.}
We can take the two limits of Example 2 and Example 3 together. 
Now the service processes will consist of weights occurring
at times of a Poisson process, with each weight having
$\Exp(1)$ distribution. In this case we get
\[
\tf(y)=\
\sup_{0<r<1} r^2\left[y-\frac1{1-r}\right].
\]
In this final case the translation into closed form is more reasonable; 
one only needs to solve a quadratic equation, 
to obtain
\[
\tf(y)=
\begin{cases}
0,& y\leq1\\
\frac1{8y^2}\Big[
8y^3+20y^2+7y+1 -\sqrt{8y+1}\left(8y^2+3y+1\right)
\Big],& y>1.
\end{cases}
\]

Taking certain other limits or special cases in these examples recovers previously known results. 
Taking $q=1-\beta$ in (\ref{generalcase}) gives the case of 
$\Geom(\beta)$ weights and leads to (\ref{Geomcase}), and further
taking $\beta\to 0$ and multiplying by $\beta$ gives the case of $\Exp(1)$ weights
and leads to (\ref{Expcase}). On the other hand, taking $\beta\to1$ in 
(\ref{generalcase}) gives Bernoulli$(q)$ weights and leads to (\ref{Bercase}).

One could also consider the continuous model of Examples 3 and 4 in the 
case where each weight has value 1,
so that the service processes are simply Poisson processes. 
To do this we take $\beta\to 1$ in Example 3, to obtain simply
\[
\tf(y)=
\Big(\big[\sqrt y-1\big]_+\Big)^2.
\]
%

Finally, by taking an appropriate continuum limit in any of these cases, one can arrive at
the Brownian first-passage percolation model which 
has been quite widely studied (see for example \cite{GUEs}, \cite{HMO}).

\section*{Acknowledgments}
Many thanks to the referee for helpful comments, and to 
Louigi Addario-Berry for suggesting reference \cite{Takacsbook}. 
Some of the work on this paper was carried out at the  
Institut Henri Poincar\'e in Paris during the programme
Interacting Particle Systems, Statistical Mechanics and Probability 
Theory; I am very grateful to the IHP for its hospitality and to INRIA 
for financial support during my stay there. The work was also supported
by an EPSRC Advanced Fellowship.

\bibliographystyle{usual2}
\bibliography{jbm}

\medskip

\noindent
\textsc{Department of Statistics,\\
University of Oxford,\\
1 South Parks Road,\\
Oxford OX1 3TG,\\
UK}\\
\texttt{martin@stats.ox.ac.uk}

\end{document}